\documentclass{amsart}
\usepackage{graphicx,amsfonts,amssymb,amsmath,amsthm}
\usepackage[pdftex]{hyperref}  
\hypersetup{pdfstartview={XYZ 120 670 1}, pdfpagemode=none}
\usepackage{cite}  
\usepackage[small,nohug,heads=vee]{diagrams}
\diagramstyle[labelstyle=\scriptstyle]

\theoremstyle{plain} 
\newtheorem{theorem}    {Theorem}[section] 

\newtheorem{corollary}  [theorem]{Corollary}

\theoremstyle{definition}
\newtheorem{definition} [theorem]{Definition}

\theoremstyle{remark}
\newtheorem{remark}              {Remark}

\newtheorem*{thm}    {Theorem} 

\numberwithin{equation}{section}

\def\C{\mathbb C}
\def\F{\mathbb F}

\def\Q{\mathbb Q}

\def\Z{\mathbb Z}

\begin{document}

\title{Matching densities for Galois representations}

\author{Nahid Walji}

\maketitle
\begin{abstract} Given a pair of $n$-dimensional complex Galois representations over $\Q$, we define their matching density to be the density, if it exists, of the set of places at which the traces of Frobenius of the two Galois representations are equal. We will show that the set of matching densities of such pairs of irreducible Galois representations (for all $n$) is dense in the interval $[0,1]$. Under the strong Artin conjecture, this also implies the corresponding statement for cuspidal automorphic representations.
\end{abstract}

\section{Introduction}

Let $\rho$, $\rho'$ be irreducible $n$-dimensional complex Galois representations over a number field $F$. For a finite place $v$ at which $\rho$ and $\rho'$ are unramified, we denote by ${\rm tr}\rho ({\rm Frob}_v)$ and ${\rm tr}\rho' ({\rm Frob}_v)$ their respective traces of Frobenius.
If $\rho$ is not isomorphic to $\rho'$, one knows that their traces of Frobenius can be equal for a set of places with an upper Dirichlet density of at most $1/2n ^2$.

Given $\rho, \rho'$, define $S := S (\rho, \rho')$ to be the set of places $v$ at which ${\rm tr}\rho ({\rm Frob}_v)={\rm tr}\rho' ({\rm Frob}_v)$. We then ask: How much can the size of $S$ vary for different pairs $(\rho, \rho')$?

We define:
\begin{definition}
Given a pair of $n$-dimensional complex Galois representations $\rho,\rho'$ over a number field $F$, the \textit{matching density} is the Dirichlet density, if it exists, of the set $\{v \mid {\rm tr}\rho ({\rm Frob}_v)={\rm tr}\rho' ({\rm Frob}_v) \}$. We denote this as $\delta (\rho,\rho')$.
We will specify \textit{non-trivial} matching densities if we mean to exclude the case where $\rho$ is isomorphic to $\rho'$.
\end{definition}

If we restrict to considering the set of non-trivial matching densities for abelian Galois representations, over any number field $F$, we note that this set is nowhere dense, being of the form $\{1/n \mid n \geq 2\}$.

There also exist examples that give values of $7/8$~\cite{Se77} and $3/4,3/5,17/32$~\cite{WaTAMS}.
In any case, one knows that for irreducible 2-dimensional complex Galois representations over any number field, the size of non-trivial matching densities is bounded above by 7/8. 
In the $n$-dimensional case, examples of Serre provide the highest possible non-trivial matching densities of $1- 1/2n ^2$~\cite{Se81}.

We will show:
\begin{theorem} \label{t1}
The set of matching densities of all pairs of irreducible $n$-dimensional complex Galois representations over $\Q$ for all $n \geq 1$ is dense in $[0,1]$.
\end{theorem}

\begin{remark}
We actually do not require all $n$ for the Theorem to hold. 
For example, as will be seen in the proof, one can restrict to square-free $n$, as well as impose any desired lower bound on $n$.
\end{remark}

Our proof of this theorem is constructive. We make use of the results of Serre~\cite{Se72} and Mazur~\cite{Ma78} on the images of $\ell$-adic representations of semistable elliptic curves, as well as the work of Iwaniec~\cite{Iw78} (via Howe--Joshi~\cite{HJ12}) on almost primes represented by quadratic polynomials, in order to construct sufficiently rich families of Galois representations.

Along the way, we will also prove a theorem about the lacunarity of $L$-functions associated to Galois representations.
Let $\delta_0 (\rho)$ denote, for an $n$-dimensional Galois representation $\rho$ over a number field, the Dirichlet density, if it exists, of places at which the associated trace of Frobenius of $\rho$ is zero. Serre's examples~\cite{Se77} provide an irreducible $n$-dimensional Galois representation $\rho$, for each positive integer $n$, such that $\delta_0 (\rho) = 1- 1/n ^2$. There are also upper bounds of $\delta_0 (\rho) \leq 1 - 1/n ^2$, for all irreducible $n$-dimensional Galois representations $\rho$~\cite{Se81}.

In this paper we prove that 
\begin{theorem} \label{t2}
The set
\begin{align*}
\{\delta_0 (\rho) \mid \text{$n$-dimensional complex Galois representations $\rho$ over $\Q$, for any }n \geq 1\}
\end{align*}
is dense in $[0,1]$.
\end{theorem}

We also consider the automorphic analogue of such questions above, often making use of the strong Artin conjecture to obtain automorphic examples from complex Galois representations. The analogous definition is:

\begin{definition}
For a pair $\pi,\pi'$ of automorphic representations for GL(n) over a number field $F$, denote their Hecke eigenvalues at a place $v$ of $F$ as $a_v(\pi), a_v(\pi')$, respectively. Then the \textit{matching density} is the Dirichlet density, if it exists, of the set $\{v \mid a_v(\pi) = a_v(\pi')\}$. We denote this as $\delta (\pi,\pi')$.
\end{definition}

For {\rm GL}(1) over any fixed number field, the set of non-trivial matching densities is $\{1/n \mid n \geq 2\}\cup \{0\}$.
In the case of {\rm GL}(2), automorphic examples arising from elliptic curves and modular forms of integral weight $\geq 2$ provide matching densities of $0,1/4,1/2,$ and $3/4$. For weight 1 cusp forms, known cases of the strong Artin conjecture give examples with values of $7/8$~\cite{Se77} and $3/4,3/5,17/32$~\cite{WaTAMS}.
For cuspidal automorphic representations for GL(2) over any number field, the size of non-trivial matching densities is bounded above by 7/8, due to Ramakrishnan's refinement of strong multiplicity one~\cite{Ra94} in this setting.

A corollary to Theorem~\ref{t1} then is:
\begin{corollary} \label{c1}
Assume the strong Artin conjecture for $n$-dimensional representations over $\Q$, for all $n$. Then the set of matching densities of pairs of cuspidal automorphic representations for ${\rm GL}(n)$ over $\Q$, for all $n \geq 1$, is dense in $[0,1]$.
\end{corollary}

We also consider the lacunarity of automorphic $L$-functions. Serre's examples correspond, via the nilpotent case of the strong Artin conjecture (due to~\cite{AC89}) to cuspidal automorphic representations $\pi_n$ for {\rm GL}(n), for each $n \geq 1$, such that $\delta_0 (\pi_n) = 1 - 1/n ^2 $. There are also upper bounds, conditional on the Ramanujan conjecture, of $\delta_0 (\pi) \leq 1 - 1/n ^2$, for all cuspidal automorphic representations $\pi$ for {\rm GL}(n)~\cite{WaLac}.
A corollary to Theorem~\ref{t2} is: 
\begin{corollary}\label{c2}
Under the strong Artin conjecture, the set
\begin{align*}
\{\delta_0 (\pi) \mid \pi \text{ a cuspidal automorphic representation for }{\rm GL}(n)/\Q, \text{ for any }n \geq 1\}
\end{align*}
is dense in $[0,1]$.
\end{corollary}

In terms of the structure of this paper, in the next section, we will provide some relevant background and construct two known examples of matching densities; in the final section, we will prove Theorems~\ref{t1} and~\ref{t2}.\\

\section{Background}

In this section, we outline the construction of two cases of matching densities that were known earlier so as to provide some context. The first case is actually a family of examples due to Serre; the second is a tetrahedral example arising in~\cite{WaTAMS}, and is one of the few examples that does not use twisting or Galois conjugacy. We also outline a conjecture of Ramakrishnan on the highest non-trivial matching density for cuspidal automorphic representations for {\rm GL}(n) over any number field.

First, we recall:
\subsection{The strong Artin conjecture}  

Consider a complex representation of the absolute Galois group over $\Q$, i.e., a continuous homomorphism 
\begin{align*}
\rho: {\rm Gal}(\overline{\Q} / \Q) \rightarrow {\rm GL}_n(\C).
\end{align*}
Due to the continuity of $\rho$ and the nature of the topologies of ${\rm Gal}(\overline{\Q} / \Q)$ and ${\rm GL}_n(\C)$, one knows that the image is always finite and the representation factors through a finite extension of ${\rm Gal}(K / \Q)$ for some $K$.

The \textit{strong Artin conjecture} implies that there exists an automorphic representation $\pi$ for GL(n) whose $L$-function is equal to that of $\rho$, and which is cuspidal if $\rho$ is irreducible. Known cases of this conjecture include the case when $n = 1$ (essentially due to the work of Hecke and Maa\ss), $n = 2$ for all representations other than those that are even with icosahedral projective image~\cite{La80, Tu81, KW109, KW209, Ki09}, and for general $n$ in the case when the image is isomorphic to a nilpotent group~\cite{AC89}.

\subsection{Examples}

We now delineate the construction of two examples of matching densities already known, the first being an example of Serre that provides matching densities of the form $1 - 1/2k ^2$ and the second being a matching density of $17/32$ from~\cite{WaTAMS}.\\

Given any positive integer $k$, there exists a representation 
\begin{align*}
\rho: {\rm Gal}(\bar{\Q} / \Q) \rightarrow {\rm GL}_k(\C),
\end{align*}
where the Dirichlet density $\delta$ is
\begin{align*}
\delta \left(\{p \mid {\rm tr}\rho ({\rm Frob}_p)= 0\} \right) = 1 - \frac{1}{k ^2}
\end{align*}
and the representation factors through a Galois group ${\rm Gal}(E / \Q)$ that is a central extension of a Galois group isomorphic to $\Z/k\Z \times \Z/k\Z:=H$. Now a central extension $G$ of $H$ satisfies the following short exact sequence 
\begin{align*}
1 \rightarrow C \rightarrow G \rightarrow H \rightarrow 1,
\end{align*}
where $C \leq Z(G)$. So $G$ has a central series $G \rhd C \rhd 1$ and thus is nilpotent.

Twist the representation $\rho$ with a suitable quadratic character $\chi$ so as to obtain a pair of complex Galois representations $(\rho, \rho \otimes \chi)$ with the desired matching density. The Galois groups that both these extensions factor through are nilpotent and so we apply the strong Artin conjecture for nilpotent groups to also obtain a pair of (cuspidal) automorphic representations with a matching density of $1-1/2k ^2$.\\

We now briefly outline an example with matching density $17/32$ and refer the reader to~\cite{WaTAMS} for further details. This will differ from the one above in that it does not rely on twisting the first representation with a character in order to obtain the second representation.

The binary tetrahedral group $\widetilde{A_4}$  is a central extension of the tetrahedral group $A_4$ by $\Z/2\Z$. This has three different irreducible 2-dimensional complex representations (up to isomorphism); exactly one of these, which we will denote as $\rho$, has traces that are all integer-valued.

Consider the following structure of number fields
 
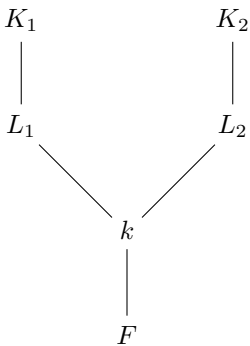
\begin{diagram}
 K_1&&&&K_2 \\
\dLine &&&& \dLine \\
 L_1&&&&L_2 \\
 &\rdLine && \ldLine & \\
& &k&& \\
& &\dLine&& \\
&& F&&
 \end{diagram}
where ${\rm Gal}(K_i / F) \simeq \widetilde{A_4}$ $(i = 1,2)$ and
 $K_1 \neq K_2$, $L_1 \neq L_2$. 
Let $\rho_1$ and $\rho_2$ be complex representations of the Galois groups ${\rm Gal}(K_1 / F)$ and ${\rm Gal}(K_2 / F)$, respectively, that correspond to $\rho$.

Galois theory and the Chebotarev density theorem show that the traces of the images of Frobenius under the complex Galois representations $\rho_1$ and $\rho_2$ are equal exactly at a set of primes of density 17/32. Applying the tetrahedral case of the strong Artin conjecture (due to Langlands~\cite{La80}) to $\rho_1$ and $\rho_2$, we obtain a pair of automorphic representations with the desired matching density.

\subsection{Conjectural upper bounds for cuspidal automorphic representations}

In general, one does expect that the examples of Serre provide the highest possible non-trivial matching densities for cuspidal automorphic representations - this is a conjecture of Ramakrishnan~\cite{Ra94b}. We remark that this does hold for $\ell$-adic Galois representations, by Rajan ~\cite{Raj98}.
Now it is known that the Ramanujan conjecture for a pair of cuspidal automorphic representations for GL(n) (if we fix any $n$) implies Ramakrishnan's conjecture for that pair as follows:

Let us assume that the distinct cuspidal automorphic representations $\pi, \pi'$ for GL(n) over a number field $F$ both satisfy the Ramanujan conjecture. We define $T = T (\pi,\pi')$ to be the set of infinite places of $F$ as well as the finite places where neither $\pi$ nor $\pi'$ are ramified. We represent their Hecke eigenvalues at a finite place $v$ of $F$ by $a_v(\pi), a_v(\pi')$, respectively.

Now consider the series 
\begin{align*}
\sum_{v \not \in T} \frac{|a_v(\pi) - a_v(\pi')|^2}{Nv^s},
\end{align*}
where the sum is over all the finite places $v$ of $F$ at which both $\pi$ and $\pi'$ are unramified. The Ramanujan conjecture implies that this series converges absolutely for ${\rm Re}(s)> 1$. For $s \in (1,2)$, we write the following inequality:
\begin{align*}
\sum_{v \not \in T} \frac{|a_v(\pi) - a_v(\pi')|^2}{Nv^s} \leq \sum_{v \in S}\frac{(2n ^2)}{Nv^s},
\end{align*}
where $S$ is the subset of places $v$ (lying in the complement of $T$) at which $a_v(\pi) \neq a_v(\pi')$. 
Now recall that the lower Dirichlet density of a set of places $P$ is defined to be
\begin{align*}
\lim_{s \rightarrow 1 ^+}{\rm inf} \frac{\sum_{v \in P}Nv^{-s}}{\sum_{v}Nv^{-s}}.
\end{align*}
We divide both sides by $\log \left(1/(s-1)\right)$ and take the limit inferior as (real) $s \rightarrow 1^+$ to obtain 
\begin{align*}
2 \leq (2n)^2 \cdot \underline{\delta}(S),
\end{align*}
where $\underline{\delta}(\cdot)$ denotes the lower Dirichlet density. Note that the value on the left-hand side of the inequality arises from the fact that, for any two cuspidal automorphic representations $\pi_1$ and $\pi_2$ for GL(n),
\begin{align*}
\sum_{v \not \in T}\frac{a_v(\pi_1) \overline{a_v(\pi_2)} }{Nv^s}= \log L^T(s, \pi_1 \times \pi_2) + O\left(1\right)
\end{align*}
under the Ramanujan conjecture, as real $s \rightarrow 1^+$. Thus
\begin{align*}
\frac{1}{2n ^2} \leq \underline{\delta}(S).\\
\end{align*}

\section{Proof of Theorems~\ref{t1} and~\ref{t2}}

We begin by constructing families of complex Galois representations that will be of use to us in the proof below.

Let $p$ be any prime greater than 7. Then one knows, due to Mazur~\cite{Ma78} and Serre~\cite{Se72},  that any semistable elliptic curve $E /\Q$ with good reduction at $p$ will have
\begin{align*}
{\rm Gal}(\Q (E[p]) / \Q) \simeq {\rm GL}_2(\F_p).
\end{align*}
The primes ramified in the extension $\Q (E[p])$ will be exactly those dividing $p N$, where $N$ is the conductor of the elliptic curve.

We construct an infinite sequence of pairwise \textit{linearly disjoint} Galois extensions (i.e., any two extensions have intersection equal to $\Q$). We achieve this by ensuring that the sets of ramified primes associated to each extension are pairwise disjoint. Therefore we will require that the conductors of the elliptic curves are all pairwise coprime. 

Howe and Joshi~\cite{HJ12} show that there exist infinitely many elliptic curves whose conductor is either a prime or a product of two primes. We want to bound the size of these primes from below. Their work crucially relies on the following theorem of Iwaniec:

\begin{thm}[Iwaniec~\cite{Iw78}]
Given any irreducible quadratic polynomial $f(x)$, there exist infinitely many positive integers $n$ such that $f(n)$ is either a prime or a product of two primes.
\end{thm}

If there are infinitely many $n$ such that $f(n)$ is prime, then this is sufficient for our purposes. So let us say that this does not hold, in which case there are infinitely many positive integers $n$ such that $f(n)$ is a product of two primes $pq$.

We want the existence of infinitely many $n$ such that both $p$ and $q$ are bounded from below. Fix some positive integer $T$ as our desired lower bound. Define the polynomial $F(n):= f (An + B)$, where $A = T!$ and some suitably chosen integer $B$ (to be specified later). Note that since $f$ is primitive, so is $F$. Writing $f(n)$ as $a n ^2 + b n + c$, we have 
\begin{align*}
F(n) = \left(a(T!n)^2 + 2aBT!n + b (T!n) \right) + \left(a B ^2 + bB + c \right).
\end{align*}
The first of the two terms on the right-hand side is divisible by all primes less than $T$, and we can choose $B$ such that the second term on the right-hand side is not divisible by any primes less than $T$; therefore, we will have that $F(n)$ will not have any prime divisors less than $T$ either.
Applying Iwaniec's theorem to $F$, we obtain infinitely many $n$ such that $F(n)$ has two prime divisors, both of which are greater than $T$.

Denote the $j$th prime as $p_j$. Fix some $\epsilon > 0$ and a point $c \in (0,1)$. We will construct a pair $\rho,\rho'$ of irreducible $n$-dimensional Galois representations over $\Q$, for some $n$, such that $\delta (\rho, \rho')$ is within $\epsilon$ of $c$. 
Pick a prime $p = p_k > 1/\epsilon$. Choose $m$ such that 
\begin{align*}
\prod_{j = k }^{ k + m}\frac{p_j-1}{p_j}
\end{align*}
is within $\epsilon$ of $c$. 

We fix suitable values of $T > p _{k + m}$ so as to obtain a sequence of $(m + 1)$ semistable elliptic curves $E_0, \dots, E_m$ with conductors $N_0, \dots, N_m$ (respectively) such that: they are a product of at most two primes, their prime divisors are all greater that $p_{k + m}$, and they are pairwise coprime. 
So we construct extensions $L_j := \Q (E_j [p_{k + j}])$ over $\Q$ which will have Galois groups ${\rm GL}_2 (\F_{p_{k + j}})$ for each $j = 0, \dots, m$. The primes ramified in the extension $L_j$ will be exactly those dividing $p_{k + j}N_j$. Since, by construction, the values $p_{k + j}N_j$ are pairwise coprime, we know that the extensions $L_j$ will have pairwise intersection $\Q$ (for all $j$).

Now ${\rm GL}_2(\F_p)$ has a $p$-dimensional irreducible representation $\rho_p$ whose associated character is equal to zero for $(1/p)|{\rm GL}_2(\F_p)|$ of the group elements. The product representation $\rho_p \times \rho_{p'}$ is an irreducible representation whose character is equal to zero on $(\frac{p-1}{p}\cdot \frac{p'-1}{p'})$ of the group elements.

We construct complex representations associated to these Galois extensions. Our extensions, being linearly disjoint, satisfy $L_{k_p}  L_{k_{p'}} \simeq {\rm GL}_2(\F_p) \times {\rm GL}_2(\F_{p'})$ when $p \neq p'$.
Thus we obtain an irreducible complex Galois representation where the proportion $w$ of Galois group elements at which the character of the representation is zero is of the form  
\begin{align*}
w(P_{k,k+m}) = \prod_{p \in P_{k,k+m}} \frac{p-1}{p},
\end{align*}
where $P_{k,k + m}$ is the set of primes $\{p_k, \dots, p_{k + m}\}$.

We note that the set $\{w(P_{k,k + m}) \mid p_k > 1/\epsilon, \ m \geq 1 \}$ is dense in $[0,1]$.
Indeed, the sequence $w(\{p_k\}), w (\{p_k, p_{k + 1}\}), w (\{p_k, p_{k+1}, p_{k+2}\}), \dots$ tends to zero and that the largest gap between any two elements of this sequence is $1/p_k$.
This proves Theorem~\ref{t2} and Corollary~\ref{c2}.

Continuing with the proof of Theorem~\ref{t1}, consider the set $S$ of all possible matching densities associated to pairs of $n$-dimensional Galois representations over $\Q$ ($n$ can vary, but is the same for each element of the pair). 

Note that the $w(P)$'s are limit points of this set. Indeed, given a complex Galois representation $\rho$ where the associated character is zero for a proportion of $w$ of the elements, if we twist this representation by a Dirichlet character $\chi$ of order $d$ (whose associated field extension is linearly disjoint to that of $\rho$), we then see that $\rho$ and $\rho \otimes \chi$ have equal traces of Frobenius exactly at a set of primes of density $w + \frac{1-w}{d}$. Letting $d$ tend to infinity, we see that $w$ is a limit point in $S$.
Thus given a point $c \in [0,1]$ and some fixed $\epsilon > 0$, we can construct a sequence that has an element within $\epsilon$ of $c$. An example of such a sequence would be the one defined above, where $k$ is chosen such that $1/p_k < \epsilon$.
This proves Theorem~\ref{t1} and Corollary~\ref{c1}.

\begin{remark}
One might ask if this method can be applied using nilpotent groups, rather than the family ${\rm GL}_2(\F_p)$, so that we might make use of the work of Arthur--Clozel~\cite{AC89} on the strong Artin conjecture for nilpotent groups, thus obtaining an unconditional result for automorphic representations.

The issue is that irreducible representations of nilpotent groups are all monomial. Thus the character formula implies that the irreducible characters of nilpotent groups vanish on at least half the group elements. On the other hand, a key property of the family ${\rm GL}_2(\F_p)$ is that it has irreducible characters that only vanish on arbitrarily small proportions of group elements. Thus using nilpotent groups is not compatible with the method above.\\
\end{remark}

\subsection{Acknowledgements:} The author would like to thank Sara Arias-de-Reyna and Dimitar Jetchev for some useful discussions, as well as Farrell Brumley and Sug Woo Shin for their helpful comments on an earlier draft of this paper.


\begin{thebibliography}{0000}
\bibitem[AC89]{AC89}
Arthur, J. and Clozel, L., Simple algebras, base change, and the advanced theory of the trace formula, Annals of Mathematics Studies, vol. 120, Princeton University Press, Princeton, NJ, 1989.

\bibitem[HJ12]{HJ12}
Howe, S., Joshi, K., Elliptic curves of almost-prime conductor, preprint (2012). 

\bibitem[Iw78]{Iw78}
Iwaniec, H.. Almost-primes represented by quadratic polynomials. Invent. Math., 47(2):171-188, 1978.

\bibitem[KWa09]{KW109}
C.~Khare and J.-P.~Wintenberger, {Serre's modularity
  conjecture. {I}}, Invent. Math. \textbf{178}, 485--504 (2009).

\bibitem[KWb09]{KW209}
\bysame, {Serre's modularity conjecture. {II}}, Invent. Math. \textbf{178}, 505--586  (2009).

\bibitem[Ki09]{Ki09}
M.~Kisin, {Modularity of 2-adic {B}arsotti-{T}ate representations},
  Invent. Math. \textbf{178}, 587--634 (2009). 

\bibitem[La80]{La80}
R.P.~Langlands, Base change for GL(2), Annals of Mathematics Studies, vol. 96, Princeton University Press, Princeton, N.J., 1980.

\bibitem[Ma78]{Ma78}
Mazur, B., Rational isogenies of prime degree, Invent. Math. 44 (1978) 129-162.

\bibitem[Raj98]{Raj98}
Rajan, C. S. \emph{On strong multiplicity one for $l$-adic representations}. Internat. Math. Res. Notices 1998, no. 3, 161--172.

\bibitem[Ram94a]{Ra94b}
D. Ramakrishnan, \emph{Pure motives and automorphic forms}, Motives
  ({S}eattle, {WA}, 1991), Proc. Sympos. Pure Math., vol.~55, Amer. Math. Soc.,
  Providence, RI, 1994, pp.~411--446. \MR{1265561 (94m:11134)}

\bibitem[Ram94b]{Ra94}
\bysame, \emph{A refinement of the strong multiplicity one theorem for {${\rm
  GL}(2)$}. {A}ppendix to: ``{$l$}-adic representations associated to modular
  forms over imaginary quadratic fields. {II}'' [{I}nvent.\ {M}ath.\ {\bf 116}
  (1994), no.\ 1-3, 619--643; {MR}1253207 (95h:11050a)] by {R}. {T}aylor},
  Invent. Math. \textbf{116} (1994), no.~1-3, 645--649. \MR{1253208
  (95h:11050b)}

\bibitem[Se72]{Se72}
Serre, J.-P.. Propri\'et\'es galoisiennes des points d'ordre fini des courbes elliptiques. Invent. Math., 15:259-331, 1972.

\bibitem[Se77]{Se77}
\bysame, Modular forms of weight one and Galois representations, Algebraic number fields: L-functions and Galois properties (Proc. Sympos., Univ. Durham, Durham, 1975), pp. 193-268, Academic Press, London, 1977.

\bibitem[Se81]{Se81}
\bysame, \emph{Quelques applications du th\'eor\`eme de densit\'e de
  {C}hebotarev}, Inst. Hautes \'Etudes Sci. Publ. Math. (1981), no.~54,
  323--401. \MR{644559 (83k:12011)}

\bibitem[Tu81]{Tu81}
J.~Tunnell, {Artin's conjecture for representations of octahedral
  type}, Bull. Amer. Math. Soc. (N.S.) \textbf{5}, 173--175 (1981).

\bibitem[Wa14a]{WaTAMS}
N. Walji, \emph{Further refinement of strong multiplicity one for {${\rm
  GL}(2)$}}, {T}ransactions of the AMS, Volume 366, Number 9, Pages 4987-5007 (2014).

\bibitem[Wa14b]{WaLac}
\bysame, \emph{On the occurence of Hecke eigenvalues and a lacunarity question of Serre}, preprint (2014).
\end{thebibliography}
\end{document}